\documentclass{commat}

\newcommand{\Z}{\mathbb{Z}}
\newcommand{\R}{\mathbb{R}}
\newcommand{\F}{\mathbb{F}}
\newcommand{\G}{\mathbb{G}}

\title{%
    Existence and uniqueness theorems for functional equations 
    }

\author{%
    Tamás Glavosits, Zsolt Karácsony
    }

\affiliation{
    \address{Tamás Glavosits --
    University of Miskolc, Department of Applied Mathematics, Miskolc-Egyetemváros, Hungary, H-3515
        }
    \email{%
    tamas.glavosits@uni-miskolc.hu
    }
    \address{Zsolt Karácsony --
    University of Miskolc, Department of Applied Mathematics, Miskolc-Egyetemváros, Hungary, H-3515
        }
    \email{%
    zsolt.karacsony@uni-miskolc.hu
    }
    }

\abstract{%
   In this paper we give simple extension and uniqueness theorems for restricted additive and logarithmic functional equations.
    }

\dedication{This article is dedicated to the memory of János Aczél (1924-2020).}



\VOLUME{32}
\NUMBER{1}
\YEAR{2024}
\firstpage{93}
\DOI{https://doi.org/10.46298/cm.10830}

\begin{paper}
 
\section{Introduction}
The main purpose of this article is to show that if $X=X(+,\leqslant)$ is an Archimedean ordered dense Abelian group, $Y$ is an Abelian group $\varepsilon\in X_{+}:=\left\{x\in X|x>0\right\}$ and  $f:\left]-2\varepsilon,2\varepsilon\right[\subseteq X\to Y$ is a function such that 
\begin{equation}\label{Equ:Add}
f(x+y)=f(x)+f(y)\qquad(x,y\in\left]-\varepsilon,\varepsilon\right[),
\end{equation}
then there uniquely exists an additive function $a:X \to Y$ such that 
\[
f(x)=a(x)\qquad x\in\left]-2\varepsilon,2\varepsilon\right[.
\]
Analogue Theorems concerning logarithmic functions are proven as well.

 Let $D\subseteq\R^{2}$ be a fixed set and define the sets $D_{x}$, $D_{y}$, $D_{x+y}$ by 
\[
\begin{aligned}
D_{x}&:=\left\{u\in\R\mid \exists(v\in\R):(u,v)\in D\right\},\\
D_{y}&:=\left\{v\in\R\mid \exists(u\in\R):(u,v)\in D\right\},\\
D_{x+y}&:=\left\{z\in\R\mid \exists((u,v)\in D):z=x+y\right\}.
\end{aligned}
\]
If the function $f:D_{x}\cup D_{y}\cup D_{x+y}\to\R$ satisfies the functional equation 
\begin{equation}
f(x+y)=f(x)+f(y)\qquad((x,y)\in D)\label{Equ:RestAdd},
\end{equation}
then the function $f$ is said to be additive on the set $D$ and the equation (\ref{Equ:RestAdd}) is said to be restricted additive functional equation. The restricted additive functional equations have previously been studied by many researchers.
In the book \cite{AL 09}  Part IV. Geometry,
Section Extension of Functional Equations p. 447--460 the authors cite numerous papers that investigate the cases when there exists an additive function $F:\R\to\R$, that is, $F(x+y)=F(x)+F(y)$ for all $x,y\in\R$ such that the function $F$ extends the function $f$, that is,  $F(x)=f(x)$ for all $x\in\mathcal{D}_{f}$ where $\mathcal{D}_{f}$ denotes the domain of the function $f$. An incomplete list of such papers is given below:
\begin{itemize}
\item In the paper \cite{AE 65} $D=(D_{+}\cup\{0\})^{2}$. 

\item In the book \cite{A 66} the first appearance of the concept of quasi-extension can be found. An additive function $a$ is said to be quasi extension of the function $f$ if $f$ is additive on a set $D\subseteq\R^{2}$ and there exist constants $c_{1}$, $c_{2}\in\R$ such that $f(u)=a(u)+c_{1}$ for all $u\in D_{x}$; $f(v)=a(v)+c_{2}$ for all $v\in D_{y}$ and $f(z)=a(z)+c_{1}+c_{2}$ for all $z\in D_{x+y}$. For example, if the function $f:\left]0,1\right[\cup\left]1,2\right[$ is defined by $f(x):=0$  whenever $x\in\left]0,1\right[$; $f(x):=1$  whenever $x\in\left]1,3\right[$, then it is easy to see that function $f$ is additive on the set $D:=\left]0,1\right[\times\left]1,2\right[$. Although $f$ has no additive extension to $\R^{2}$ exists, to identically zero function is an additive quasi extension of the function $f$ from $D$ to $\R^{2}$.   

\item In the paper \cite{DL 67} the cases are investigated  when $D=\R_{+}^{2}$ and $D$ is an open interval of the real line containing the origin. In this paper the notations $D_{x}$, $D_{y}$, $D_{x+y}$ has appeared first.

\item In \cite{Sz 72} the author generalizes the above result that $D\subseteq \R^{2}$ is an arbitrary open set, $D_{0}=D_{x}\cup D_{y}\cup D_{x+y}$, $f:D_{0}\to \R$ is a function such that $f(x+y)=f(x)+f(y)$ for all $(x,y)\in D$.

\item In \cite{Ri 76} a simple extension theorem  can be found for Pexider additive functional equation where the additivity is fulfilled in a nonempty connected open set of the real line.

\item In the article \cite{A 83} $D=H(I)$ where $I$ is a nonempty open interval of the real line and the set $H(I)$ is defined by 
\[
H(I):=\left\{(x,y)\in\R^{2}\mid x,y,x+y\in I\right\}.
\]
The set $H(I)$ is a hexagon, sometimes a triangle or the empty set.

\item In the book \cite{K 08} $D\subseteq\R^{N}$ is a nonempty connected open set. The extension is brought back to the theory of convex functions, but in this book the author does not consider the restricted Pexider additive functional equations. 

\item In the article \cite{RB 87} an extension theorem can be found for restricted Pexider additive functional equations where $D\subseteq\R^{N}$ is a nonempty connected open set.

\item In the book \cite{AD 89} several functional equations are considered in more general abstract algebraic settings.
\end{itemize}

 Below some necessary concepts and notations are collected:

 Let $X=X(+,\leqslant)$ be an ordered group. The absolute value of an element $x\in X$ is defined by 
\[
|x|:=\max\{x,-x\}.
\]

 Let $X=X(+,\leqslant)$ be an ordered group. If $x\in X$ and $n\in\Z$ where $\Z$ denotes the ring of integers, then we can define the element $nx\in X$ by 
\[
nx:=
\left\{
\begin{array}{cl}
x+\dots+x,&\hbox{ if $n>0$};\\
0,&\hbox{ if $n=0$};\\
(-x)+\dots+(-x),&\hbox{ if $n<0$}.
\end{array}
\right.
\]

 An ordered group $X=X(+,\leqslant)$ is said to be Archimedean ordered, if for every two elements $x, y\in X_{+}:=\left\{z\in X|z>0\right\}$, there exists $n\in\Z_{+}:=\left\{1,2,\dots\right\}$, such that $y<nx$.

 An ordered group $X=X(+,\leqslant)$ is said to be dense, if
\[
\left]a,b\right[:=\left\{x\in X\mid a<x\text{ and }x<b\right\}\neq\emptyset
\]
for all $a,b\in X$ with $a<b$.

 An ordered field $\F=\F(+,\cdot,\leqslant)$ is said to be Archimedean ordered, if the ordered group $\F=\F(+,\leqslant)$ is Archimedean ordered. (We use the concept of the field including the commutativity of the  operation '$\cdot$'.)

 A homomorphism $a:X(+)\to Y(+)$, that is, a function $a:X\to Y$ with
\[
a(x+y)=a(x)+a(y)\qquad(x,y\in X)
\]
is said to be an additive function.

 A homomorphism $l:X(\cdot)\to Y(+)$, that is, a function $l:X\to Y$ with
\[
l(xy)=l(x)+l(y)\qquad(x,y\in X)
\] 
is said to be a logarithmic function.

 In the rest of this article we use four properties of the open intervals in the appropriately ordered structure \cite{GK 02}:
\begin{enumerate}
\item If $\G(+,\leqslant)$ is an ordered group, then the open intervals are translation invariant, that is, 
\[
\gamma+\left]\alpha,\beta\right[=\left]\gamma+\alpha,\gamma+\beta\right[\qquad(\gamma\in\G).
\]
\item If $\G(+,\leqslant)$ is an ordered dense Abelian group,  $\alpha$, $\beta$, $\gamma$, $\delta\in\G$ such that $\alpha<\beta$ and $\gamma<\delta$, then 
\[
\left]\alpha,\beta\right[+\left]\gamma,\delta\right[ =\left]\alpha+\gamma,\beta+\delta\right[.
\]
\item If $\F(+,\cdot,\leqslant)$ is an ordered field, then the open intervals are homothety invariant, that is, if  $\alpha$, $\beta$, $\gamma\in\F$ such that $\alpha<\beta$ and $\gamma>0$, then 
\[
\gamma\cdot\left]\alpha,\beta\right[=\left]\gamma\alpha,\gamma\beta\right[.
\]
\item If $\F(+,\cdot,\leqslant)$ is an ordered field, $\alpha$, $\beta$, $\gamma$, $\delta\in\F$ such that $0<\alpha<\beta$ and $0<\gamma<\delta$, then 
\[
\left]\alpha,\beta\right[\cdot\left]\gamma,\delta\right[ =\left]\alpha\gamma,\beta\delta\right[.
\]
\end{enumerate}
Property (2) can be easily deduced from property (1), although in  \cite{GK 01} can be found an example of dense Abelian semigroup which has property (2) without property (1). 

 Similarly, property (4) can be easily deduced from property (3), although in  \cite{GK 01} can also be found an example for dense Abelian semigroup which has property (4) without property (3).

 Our paper  is structured as follows: In section 2 we consider the additive and multiplicative versions of Euclid's Theorem, which will be the key to our extension theorems for additive and logarithmic functions respectively. In section 3 we give extension theorems for additive and logarithmic functions. In section 4 we give uniqueness theorems for additive and logarithmic functions.

\section{Euclid's Theorem}

Euclid's  Elements \cite{Eu} is one of the most influential mathematical textbooks written more than two thousand years ago. In this textbook (book X, proposition 3.) there is an algorithm using the so-called Euclidean or remainder division to give the greatest common measure of two given commensurable magnitudes. We use the modern version of this division to give our extension theorem for restricted additive functional equations. We start with the existence and uniqueness theorem of Euclidean division. 

\begin{theorem}\label{Thm:EucDiv}
If $\G=\G(+,\leqslant)$ is an Archimedean ordered group, $x$, $y\in \G$ with $y\neq 0$. Then there uniquely exists an integer $q$ and an element $r\in \G$ such that 
\[
x=qy+r\qquad\text{where}\qquad0\leqslant r<|y|.
\]
\end{theorem}

\begin{proof}
It is easy to see that 
\[
\G=\bigcup_{z\in\Z}\left(zy+\left[0,|y|\right[\right)
\]
and the union is disjoint whence the theorem is clear.
\end{proof}

\begin{proposition}\label{Prop:Berno}
If $\F$ is an Archimedean ordered field, $x$, $y\in\F_{+}$, $x\neq 0$. Then there exists $n$, $m\in\Z_{+}$ such that $x^{m}<y<x^{n}$.
\end{proposition}

\begin{proof}
First, we investigate the case when $x>1$.
Let $h:=x-1$. Thus we have that there exists $n\in\Z_{+}$ such that 
\[
y<1+nh
\]
and
by the Bernoulli inequality \cite{B 89} we have that 
\[
1+nh\leqslant (1+h)^{n}=x^{n}.
\]
From the above two inequalities we obtain, that $y<x^{n}$. 

 It is also easy to see that 
\[
\lim_{n\to\infty}\left(\frac{1}{x}\right)^{n}=0,
\]
thus we obtain that there exists an integer $m$ such that $x^{m}<y$.
The case $0<x<1$ is similar.
\end{proof}

Now we show the multiplicative version of the Euclidean division.

\begin{theorem}\label{Thm:EucDivMult}
If $\F=\F(+,\cdot,\leqslant)$ is an Archimedean ordered field, $x$, $y\in\F_{+}$ such that $y\neq 1$. Then there uniquely exists an integer $z$ and an element $r\in \F_{+}$ such that 
\[
x=y^{z}\cdot r.
\]
Furthermore, if $1<y$, then $1<r<y$; and if $y<1$, then  $y<r<1$.
\end{theorem}

\begin{proof}
By Proposition \ref{Prop:Berno}. we have that 
\[
\F_{+}=\left\{
\begin{array}{ll}
\displaystyle\bigcup_{z\in\Z}y^{z}\left[1,y\right[&\text{ if }1<y;\\
\\
\displaystyle\bigcup_{z\in\Z}y^{z}\left[y,1\right[&\text{ if }y<1;
\end{array}
\right.
\]
and the union is disjoint, whence the proof can be easily derived.
\end{proof}

\section{Extension Theorems for additive and logarithmic functional equations}

\begin{theorem}\label{Thm:KiterjesztII1}
Let $\G(+,\leqslant)$ be an Archimedean ordered dense Abelian group, $Y(+)$ be a group, $\varepsilon\in \G_{+}$ and
 $f:\left]-2\varepsilon,2\varepsilon\right[\to Y$ be a function such that
\begin{equation}\label{Equ:fadd}
f(x+y)=f(x)+f(y)\qquad (x,y\in\left]-\varepsilon,\varepsilon\right[),
\end{equation}
then there exists an additive function $a:\G\to Y$ which extends the function $f$.
\end{theorem}

\begin{proof}
Define the function $a:\G\to Y$ by
\[
a(x):= nf(y_{0})+f(r)
\]
where $y_{0}\in \left]0,\varepsilon\right[$ is an arbitrarily fixed element and the element $x\in \G$ is of the form  $x=ny_{0}+r$ where $n\in\Z$ and $r\in \G$ such that $0\leq r<y_{0}$. This form of $x$ is unique by Theorem \ref{Thm:EucDiv}.

 We show that the function $a$ is additive. For this let $x,y\in \G$. By Theorem \ref{Thm:EucDiv}. we have that
\begin{equation}\label{Equ:Euk1}
\begin{aligned}
x&=n_{1}y_{0}+r_{1},\\
y&=n_{2}y_{0}+r_{2},
\end{aligned}
\qquad
\begin{aligned}
0&\leqslant r_{1}<y_{0}\\
0&\leqslant r_{2}<y_{0}.
\end{aligned}
\end{equation}
Then
\[
x+y=(n_{1}+n_{2})y_{0}+r_{1}+r_{2}=(n_{1}+n_{2})y_{0}+r_{3},
\]
where $r_{3}:=r_{1}+r_{2}$ thus $0\leqslant r_{3}<2y_{0}$. There are two cases:
\begin{itemize}
  \item If $0\leqslant r_{3}<y_{0}$, then
  \[
  \begin{aligned}
    a(x+y)&=(n_{1}+n_{2})f(y_{0})+f(r_{3})=\\
    \quad&(n_{1}+n_{2})f(y_{0})+f(r_{1}+r_{2})\stackrel{(\ref{Equ:fadd})}{=} \\
    \quad&[n_{1}f(y_{0})+f(r_{1})]+[n_{2}f(y_{0})+f(r_{2})]=a(x)+a(y).
  \end{aligned}
  \]
  \item If $y_{0}\leqslant r_{3}<2y_{0}$, that is, $0\leqslant r_{3}-y_{0}<y_{0}$. Then
  \begin{equation}\label{Eku:Euk21}
    x+y=(n_{1}+n_{2}+1)y_{0}+(r_{3}-y_{0}),
  \end{equation}
  and the another hand since $r_{3}-y_{0}\in\left]\varepsilon,\varepsilon\right[$ and $y_{0}\in \left]-\varepsilon,\varepsilon\right[$ thus
  \begin{equation}\label{Eku:Euk22}
  f(r_{3})=f(r_{3}-y_{0})+f(y_{0}).
  \end{equation}
  Whence we have that
  \begin{align*}
    a(x+y)
    { }&{ }\stackrel{(\ref{Eku:Euk21})}{=}(n_{1}+n_{1}+1)f(y_{0})+f(r_{3}-y_{0}) \\
    &= (n_{1}+n_{2})f(y_{0})+[f(r_{3}-y_{0})+f(y_{0})] \\
    &\stackrel{(\ref{Eku:Euk22})}{=} (n_{1}+n_{2})f(y_{0})+f(r_{3}) \\
    &= (n_{1}+n_{2})f(y_{0})+f(r_{1}+r_{2}) \\
    &\stackrel{(\ref{Equ:fadd})}{=} (n_{1}+n_{2})f(y_{0})+f(r_{1})+f(r_{2}) \\
    &= [n_{1}f(y_{0})+f(r_{1})]+[n_{1}f(y_{0})+f(r_{1})] \\
    &= a(x)+a(y).
  \end{align*}
\end{itemize}

 We show that 
\[
f(x)=a(x)\qquad (x\in\left]-\varepsilon,\varepsilon\right[).
\]
For this let $x\in \left]-\varepsilon,\varepsilon\right[$. Then by Euclidean division, we obtain that there exists a number $n\in \Z$ and an element $r\in \G$ such that
\[
x=ny_{0}+r\qquad\text{where}\qquad 0\leqslant r<y_{0}.
\]
There are three cases:
\begin{itemize}
 \item If $x\in\left[0,\varepsilon\right[$, then 
  \[
  a(x)=nf(y_{0})+f(r)\stackrel{(\ref{Equ:fadd})}{=}f(ny_{0})+f(r)\stackrel{(\ref{Equ:fadd})}{=}f(ny_{0}+r)=f(x).
  \]
  \item If $x\in\left]-\varepsilon,0\right[$ and $ny_{0}\in\left]-\varepsilon,0\right[$, then we can apply the chain of reasoning of the first case. 
  \item If $x\in\left]-\varepsilon,0\right[$ and $ny_{0}<-\varepsilon$, then 
  \[
  ny_{0}<x<(n+1)y_{0},\quad\text{and}\quad (n+1)y_{0}\in\left]-\varepsilon,0\right[,
  \]
  whence we have that
  \begin{equation}\label{Equ:Euk23}
    (n+1)f(y_{0})=f((n+1)y_{0}).
  \end{equation}

  On the other hand $(y_{0}-r)\in ]0,\varepsilon[$, $r\in\left]0,\varepsilon\right[$ thus $f(y_{0})=f(y_{0}-r)+f(r)$ but 
  $f(y_{0}-r)=-f(r-y_{0})$ whence we have that 
  \begin{equation}\label{Equ:Euk24}
  f(r)=f(r-y_{0})+f(y_{0}).
  \end{equation}
  Thus we obtain that
  \begin{align*}
  a(x)
  { }&{ }= nf(y_{0})+f(r) \\
  &\stackrel{(\ref{Equ:Euk24})}{=} nf(y_{0})+f(r-y_{0})+f(y_{0}) \\
  &= (n+1)f(y_{0})+f(r-y_{0}) \\
  &\stackrel{(\ref{Equ:Euk23})}{=} f((n+1)y_{0})+f(r-y_{0}) \\
  &\stackrel{(\ref{Equ:fadd})}{=} f(ny_{0}+r)=f(x).
  \end{align*}
\end{itemize}

 Finally we show that 
\[
f(x)=a(x)\qquad (x\in\left]-2\varepsilon,2\varepsilon\right[).
\]
For this let $x\in\left]-2\varepsilon,2\varepsilon\right[$. Then based on the relation for the sum of the intervals we get that
\[
\left]-2\varepsilon,2\varepsilon\right[=\left]-\varepsilon,\varepsilon\right[+\left]-\varepsilon,\varepsilon\right[
\]
thus there exist elements $u,v\in \left]-\varepsilon,\varepsilon\right[$ such that $x=u+v$.
Thus by the previous part of this proof we get
\[
f(x)=f(u+v)\stackrel{(\ref{Equ:fadd})}{=}f(u)+f(v)=a(u)+a(v)=a(u+v)=a(x).
\qedhere
\]
\end{proof}

\begin{theorem}\label{Thm:ExtLog}
Let $\F(+,\cdot,\leqslant)$ be an Archimedean ordered field, $Y(+)$ be a group, $\varepsilon\in\F_{+}$,
 $f:\left]\varepsilon^{-2},\varepsilon^{2}\right[\to Y$ be a function such that
\begin{equation}\label{Equ:flog}
f(xy)=f(x)+f(y)\qquad(x,y\in \left]\varepsilon^{-1},\varepsilon\right[),
\end{equation}
then there exists a logarithmic function $l:\F_{+}\to Y$ which extends the function $f$.
\end{theorem}

\begin{proof}
Define the function $l:\F_{+}\to Y$ by
\[
l(x):= nf(y_{0})+f(r)
\]
where $y_{0}\in \left]1,\varepsilon\right[$ is an arbitrarily fixed element and the element $x\in \F_{+}$ is of the form  $x=y_{0}^{n}\cdot r$ where $n\in\Z$ and $r\in \F_{+}$ such that $1\leq r<y_{0}$. This form of $x$ is unique by Theorem \ref{Thm:EucDivMult}. Similarly as in the proof of the Theorem \ref{Thm:KiterjesztII1}. it is easy to show that the function $l$ is logarithmic and extends the function $f$.
\end{proof}

\section{Uniqueness Theorem for additive and logarithmic functional equations}

\begin{theorem}\label{Thm:Unic2}
Let $\G$ be an Archimedean ordered Abelian group and $a:\G\to Y$ be an additive function. If there exist constants $\alpha$ $\beta\in\G$ with $\alpha<\beta$ and $c\in Y$  such that 
\[
a(x)=c\qquad (x\in\left]\alpha,\beta\right[),
\]
then $a(x)=0$ for all $x\in\G$.
\end{theorem}

\begin{proof}
By the translation invariant property of intervals, it is easy to see that $a(y)=d$ for all $y\in\left]0,\beta-\alpha\right[$ where $d=c-a(\alpha)$. Let $\varepsilon\in\left]0,\beta-\alpha\right[$ and $\delta\in\left]0,\varepsilon\right[$. Then 
\[
d=a(\varepsilon)=a(\varepsilon-\delta+\delta)=a(\varepsilon-\delta)+a(\delta)=d+d
\]
thus we have that $a(y)=0$ for all $y\in\left]0,\beta-\alpha\right[$.

 Let $x\in\G$ be arbitrary. By the Theorem \ref{Thm:EucDiv}. there exists an integer $z$ and an element $r\in \G$ such that $0\leqslant r<\varepsilon$ and $x=q\varepsilon+r$, whence we obtain that
\[
a(x)=a(q\varepsilon+r)=q\cdot a(\varepsilon)+a(r)=q\cdot 0+0=0
\]
which completes the proof.
\end{proof}

\begin{corollary}
Let $\G$ be an Archimedean ordered Abelian group and $a_{1}$, $a_{2}:\G\to Y$ be additive functions. If there exists a nonempty open interval $\left]\alpha,\beta\right[\subseteq X$  and a constant $c\in Y$ such that
\[
a_{1}(x)=a_{2}(x)+c\qquad (x\in\left]\alpha,\beta\right[),
\]
then $a_{1}(x)=a_{2}(x)=0$ for all $x\in \G$.
\end{corollary}

\begin{theorem}
Let $\F(+,\cdot,\leqslant)$ be an Archimedean ordered field, $Y(+)$ be a  group, $l:\F_{+}\to Y$ be a logarithmic function. If there exists a nonempty internal $\left]a,b\right[\subseteq \F_{+}$ and a constant $c\in Y$ such that
\[
l(x)=c\qquad (x\in\left]a,b\right[),
\]
then $l(x)=0$ for all $x\in\F_{+}$.
\end{theorem}

\begin{proof}
By the homothety invariant property of intervals, it is easy to see that $a(y)=d$ for all $y\in\left]1,\frac{\beta}{\alpha}\right[$ where $d=c-l(\alpha)$. Let $\varepsilon\in\left]1,\frac{\beta}{\alpha}\right[$ and $\delta\in\left]1,\varepsilon\right[$, then 
\[
d=a(\varepsilon)=a\left(\frac{\varepsilon}{\delta}\cdot\delta\right)=a\left(\frac{\varepsilon}{\delta}\right)+a(\delta)=d+d
\]
thus we have that $a(y)=0$ for all $y\in\left]1,\frac{\beta}{\alpha}\right[$.

 Let $x\in \F_{+}$ be arbitrary. By the Theorem \ref{Thm:EucDivMult}. there exists an integer $z$ and an element $r\in \F_{+}$ such that $1\leqslant r<\varepsilon$ and $x=\varepsilon^{q}\cdot r$ whence we obtain that
\[
l(x)=l(\varepsilon^{q}\cdot r)=q\cdot l(\varepsilon)+a(r)=q\cdot 0+0=0
\]
which completes the proof.
\end{proof}

\begin{corollary}
Let $\F(+,\cdot,\leqslant)$ be an Archimedean ordered group and $l_{1}$, $l_{2}:\F_{+}\to Y$ be an additive functions. If there exists a nonempty open interval $\left]\alpha,\beta\right[\subseteq \F_{+}$  and a constant $c\in Y$ such that
\[
l_{1}(x)=l_{2}(x)+c\qquad (x\in\left]\alpha,\beta\right[),
\]
then $l_{1}(x)=l_{2}(x)$ for all $x\in\F_{+}$.
\end{corollary}


\EditInfo{August 04, 2020}{March 16, 2021}{Attila Berczes}

\end{paper}